\documentclass{amsart}
\usepackage{amssymb}
\usepackage{amsmath}
\usepackage{latexsym}
\usepackage{amscd}
\usepackage{eufrak}
\usepackage{mathrsfs}
\usepackage[dutch,english]{babel}

\DeclareMathAlphabet{\mathpzc}{OT1}{pzc}{m}{it}

\newtheorem{theorem}{Theorem}[section]
\newtheorem{theorem-definition}[theorem]{Theorem-Definition}
\newtheorem{lemma-def}[theorem]{Lemma-Definition}
\newtheorem{definition-prop}[theorem]{Proposition-Definition}

\newtheorem{prop}[theorem]{Proposition}
\newtheorem{lemma}[theorem]{Lemma}
\newtheorem{cor}[theorem]{Corollary}
\newtheorem{definition}[theorem]{Definition}
\newtheorem{notation}[theorem]{Notation}

\newtheorem{conjecture}[theorem]{Conjecture}

\newenvironment{remark}{\vspace{4pt}\noindent\textbf{Remark.}}{\qed\vspace{4pt}}
\newenvironment{example}{\vspace{4pt}\noindent\textbf{Example.}}{\qed\vspace{4pt}}

\newcommand{\Spec}{\ensuremath{\mathrm{Spec}\,}}

\newcommand{\LL}{\ensuremath{\mathbb{L}}}

\newcommand{\N}{\ensuremath{\mathbb{N}}}
\newcommand{\Z}{\ensuremath{\mathbb{Z}}}
\newcommand{\Q}{\ensuremath{\mathbb{Q}}}

\newcommand{\R}{\ensuremath{\mathbb{R}}}
\newcommand{\C}{\ensuremath{\mathbb{C}}}

\newcommand{\Pro}{\ensuremath{\mathbb{P}}}
\newcommand{\A}{\ensuremath{\mathbb{A}}}

\renewcommand{\R}{\ensuremath{\mathbb{R}}}
\renewcommand{\C}{\ensuremath{\mathbb{C}}}
\renewcommand{\Pro}{\ensuremath{\mathbb{P}}}
\renewcommand{\A}{\ensuremath{\mathbb{A}}}

\numberwithin{equation}{section}

\begin{document}
\title[$P$-adic and motivic zeta functions]{An introduction to $P$-adic and motivic zeta
functions and the monodromy conjecture}
\author{Johannes Nicaise}
\address{Universit\'e Lille 1\\
Laboratoire Painlev\'e, CNRS - UMR 8524\\ Cit\'e Scientifique\\59655 Villeneuve d'Ascq C\'edex \\
France} \email{johannes.nicaise@math.univ-lille1.fr}
 \maketitle

\section{Introduction}\label{sec-intro}
Introduced by Weil, the $p$-adic zeta function associated to a
polynomial $f$ over $\Z_p$ was systematically studied by Igusa in
the non-archimedean wing of his theory of local zeta functions,
which also includes archimedean (real and complex) zeta functions
\cite{Igusa3}\cite{Igusa:intro}. The $p$-adic zeta function is a
meromorphic function on the complex plane, and contains
information about the number of solutions of the congruences
$f\equiv 0\!\!\mod p^m$ for $m>0$. Igusa formulated an intriguing
conjecture, the monodromy conjecture, stating that, if $f$ is
defined over $\Z$, the poles of its $p$-adic zeta function are
closely related to the structure of the singularities of the
complex hypersurface defined by $f$ (see Conjecture \ref{con-mon}
for the precise statement). Special cases of the conjecture have
been solved (in particular the case where $f$ is a polynomial in
two variables) but the general case remains quite mysterious.

\begin{definition}[Igusa Poincar\'e series]
Let $p$ be a prime. For any pair of integers $m,d> 0$, and any
polynomial $f\in \Z_p[x_1,\ldots,x_d]$, we denote by $N_{m}(f)$
the number of solutions of the congruence $f\equiv 0\mod p^{m}$ in
the ring $(\Z/(p^{m}))^{d}$. We denote by $Q(f;t)$ the generating
series
$$Q(f;t)=\sum_{m> 0}N_{m}(f)\cdot  t^m\quad \in \Z[[t]]$$
and we call it the Igusa Poincar\'e series associated to $f$.
\end{definition}

We will see that  the series $Q(f;t)$ is rational, i.e. belongs to
the subring $\Q(t)$ of $\Q[[t]]$. This was conjectured by Borevich
and Shafarevich and proven by Igusa. Igusa's $p$-adic zeta
function $Z(f;s)$ associated to $f$ is given by
$$Z(f;s)=1+\frac{(p^{-s}-1)Q(f;p^{-s-d})}{p^{-s}}$$

\begin{notation}
If $f\in \Z[x_1,\ldots,x_d]$, we can view $f$ as a $p$-adic
polynomial for any prime $p>0$, and we write $N_{p,m}(f)$,
$Q_{p}(f;t)$ and $Z_p(f;s)$ for the objects introduced above.
\end{notation}

We adopt the following notation: if $f$ is a polynomial in
$\Z[x_1,\ldots,x_d]$ then we denote by $V_{f}$ the closed
subscheme of $\A^d_{\Z}$ defined by $f$, and by
$|V_{f}(\mathbb{F}_p)|$ the cardinality of its set of
$\mathbb{F}_p$-rational points. In down-to-earth terms, this is
the number of solutions of the equation $f=0$ over the field
$\mathbb{F}_p$.

Denote by $v_p:\Z\rightarrow \N\cup\{\infty\}$ the $p$-adic
valuation and by $\overline{ac}_p:\Z_{\neq 0}\rightarrow
\mathbb{F}_p$ the reduced angular component map, i.e.
$\overline{ac}_p(z)$ is the reduction of $z\cdot p^{-v_p(z)}$
modulo $p$ for any $z\in \Z_{\neq 0}$. By convention,
$\overline{ac}_p(0)=0$. For any integer $m>0$, the maps $v_p$ and
$\overline{ac}_p$ induce maps
$$v_p:\Z/(p^m)\rightarrow \{0,\ldots,m-1\}\cup\{\infty\}$$ and
$\overline{ac}_p:\Z/(p^m)\rightarrow \mathbb{F}_p$ in the obvious
way, by composition with an arbitrary section $\Z/(p^m)\rightarrow
\Z$ of the projection $\Z\rightarrow \Z/(p^m)$ which sends $0\in
\Z/(p^m)$ to $0\in \Z$.

Let us look at some basic examples of $p$-adic zeta functions,
with $d=2$. For notational convenience we write $(x,y)$ instead of
$(x_1,x_2)$. Consider the polynomials $g(x,y)=y-x^2$ and
$h(x,y)=y^2-x^3$ in $\Z[x,y]$, and fix a prime $p$.

It is obvious that $N_{p,m}(g)=p^{m}$ for each $m>0$: we can
freely choose a value for $x$ in $\Z/(p^m)$, and this choice
determines $y$. Let us rephrase this a little bit to obtain a
formula which will turn out to generalize. For any integer $m>0$
and any solution $u$ of the equation $f=0$ in $(\Z/(p^m))^2$,
there exist exactly $p$ solutions of the equation $f= 0$ in
$(\Z/(p^{m+1}))^2$ which are mapped to $u$ under the projection
$(\Z/(p^{m+1}))^2\rightarrow (\Z/(p^{m}))^2$ (choose an arbitrary
lifting of the $x$-coordinate to $\Z/(p^{m+1})$~; this uniquely
determines the $y$-coordinate of the lifting of $u$). Therefore,
we get
\begin{eqnarray*}
Q_{p}(g;t)&=&\frac{|V_g(\mathbb{F}_p)|t}{1-pt}
\\
Z_p(g;s)&=&p^{-2}|V_g(\mathbb{F}_p)|(p-1)\frac{p^{-1-s}}{1-p^{-1-s}}+1-p^{-2}|V_g(\mathbb{F}_p)|
\end{eqnarray*}
We have $|V_g(\mathbb{F}_p)|=p$.


The polynomial $h(x,y)$ is more tricky. We exclude the case $p=3$.
Fix an integer $m>0$. Choose a value $a$ for $x$ in $\Z/(p^m)$.
The equation $y^3=a^2$ has a solution in $\Z/(p^m)$ iff
\begin{enumerate}
\item $2v_p(a)\geq m$, or \item  $2v_p(a)<m$, $v_p(a)$ is
divisible by $3$, and $\overline{ac}_p(a^2)$ is a cube in
$\mathbb{F}_p$.\end{enumerate} In the first case, we can take for
$y$ any element $b$ in $\Z/(p^m)$ with $3v_p(b)\geq m$. In the
second case, we put $\omega=v_p(a)/3$. For any solution
$\overline{b}$ of the equation $y^3-\overline{ac}_p(a^2)=0$ in
$\mathbb{F}_p$, there exists an element $b_0$ in $\Z/(p^m)$ such
that $\overline{ac}_p(b_0)=\overline{b}$ and such that the set of
all solutions $b$ of the equation $y^3-a^2=0$ in $\Z/(p^m)$ with
$\overline{ac}_p(b)=\overline{b}$ is given by the coset
$$b_0+p^{m-4\omega}\cdot \Z/(p^m)$$
Direct computation yields
\begin{eqnarray*}
Q_{p}(h;t)&=&\frac{|V_h(\mathbb{F}_p)|t+p(p-1)t^2+(p-1)p^6t^6-p^8t^7}{(1-pt)(1-p^7t^6)}
\\
Z_p(h;s)&=&1+\frac{(p^{-s}-1)(p^{-2}|V_h(\mathbb{F}_p)|+(p-1)p^{-s-3}+(p-1)p^{-5s-6}-p^{-6s-6})}{(1-p^{-s-1})(1-p^{-6s-5})}
\end{eqnarray*}
and we have $|V_h(\mathbb{F}_p)|=p$. We leave the case $p=3$ as an
exercise to the reader.

These examples illustrate some general phenomena. If $f$ is a
polynomial in $\Z[x_1,\ldots,x_d]$ and $p$ is a prime such that
$p$ does not divide $f$ and such that the equation $f=0$ defines a
smooth subvariety of the affine space $\A^d_{\mathbb{F}_p}$, we
have
\begin{eqnarray*}
Q_{p}(f;t)&=&\frac{|V_f(\mathbb{F}_p)|t}{1-p^{d-1}t}
\\ Z_p(f;s)&=&p^{-d}|V_f(\mathbb{F}_p)|(p-1)\frac{p^{-1-s}}{1-p^{-1-s}}+1-p^{-d}|V_f(\mathbb{F}_p)|
\end{eqnarray*}
In fact, one can show that for any $m>0$ and any solution $u$ of
the equation $f= 0$ in $(\Z/(p^m))^d$, there exist exactly
$p^{d-1}$ solutions of the equation $f= 0$ in $(\Z/(p^{m+1}))^d$
which are mapped to $u$ under the projection
$(\Z/(p^{m+1}))^d\rightarrow (\Z/(p^{m}))^d$. For $d=1$ this
statement is well-known as Hensel's Lemma. So in the smooth case,
Igusa's zeta function is easy to compute (provided you know the
number $|V_f(\mathbb{F}_p)|$, that is!) and only has the so-called
trivial poles $s\in \{1+(2\pi i/\ln p)\Z\}$, each with
multiplicity one.

If the equation $f=0$ defines a smooth subvariety in $\A^d_{\C}$,
then the equation $f=0$ also defines a smooth subvariety in
$\A^d_{\mathbb{F}_p}$, for $p\gg 0$. So we see that in this case,
the set of poles of $Z_p(f;s)$ is $\{1+(2\pi i/\ln p)\Z\}$ for
 $p\gg 0$. This is a (very) special case of Igusa's Monodromy
 Conjecture. If the closed subscheme of $\A^d_{\C}$ defined by $f$
 is not smooth (as in our second example), the behaviour of the
 $p$-adic zeta function is much harder to control. According to
 the Monodromy Conjecture, part of this behaviour can be described
 by analyzing the singularities of the complex hypersurface
 defined by $f$. A powerful method to compute the $p$-adic zeta function is taking an
 embedded resolution of singularities for
 $f$, which essentially reduces the computation to the case where $f$ is a
 monomial; see Section \ref{subsec-rat}.

A second observation is that our examples show a certain
uniformity in $p$. For fixed $f$, there exist algebraic varieties
$V_1,\ldots,V_r$ defined over $\Z$ and rational functions
$G_1,\ldots,G_r$ in $\Q(u,v)$ of the form
$$G_{i}(u,v)=\prod_{j=1}^{n_i}\frac{u^{-a_{i,j}}v^{b_{i,j}}}{1-u^{-a_{i,j}}v^{b_{i,j}}}$$
with $a_{i,j},\,b_{i,j}\in \Z_{>0}$ such that for $p\gg 0$,
$$Z_p(f;s)=p^{-d}\sum_{i=1}^{r}|V_i(\mathbb{F}_p)|G_i(p,p^{-s})$$
 This is also a general
phenomenon: as we will see, there exists a ``universal'' zeta
function $Z_{mot}(f;s)$ associated to $f$, which is built from
algebraic varieties and specializes to the $p$-adic zeta function
for $p\gg 0$ (Theorem \ref{specialize}). This zeta function is the
so-called motivic zeta function of $f$.

\subsection*{Further reading.}
 The basic references for the theory of $p$-adic
zeta functions (Sections \ref{sec-general} and \ref{sec-pmoncon})
are \cite{Igusa:intro} and \cite{DenefBour}. For more background
on the motivic side of the story (Sections \ref{sec-motint} and
\ref{sec-motzeta}) we refer to
\cite{DL3}\cite{loeser-seattle}\cite{Vemot}.

\subsection*{Acknowledgements} I would like to thank the
referee for his careful reading of the paper.
\section{Generalities on $p$-adic zeta
functions}\label{sec-general}
\subsection{Definitions}\label{subsec-def}
Throughout Sections \ref{sec-general} and \ref{sec-pmoncon}, we
fix the following notations: $K$ is a $p$-adic field (i.e. a
finite extension of $\Q_p$ for some prime $p$), $R$ its ring of
integers, $P$ the maximal ideal of $R$, and $k=R/P$ the residue
field. We fix a generator $\pi$ of $P$, and we denote the
cardinality of $k$ by $q$.

 We denote by
$v(\cdot)$ the discrete ($\pi$-adic) valuation on $K$, and we put
$|z|_K=q^{-v(z)}$ for $z\in K^*$. We denote by $ac(z)$ the angular
component $z\cdot\pi^{-v(z)}$ of $z\in K^*$; it is an element of
the group of units $R^{\ast}$. We put $|0|_K=0$ and $ac(0)=0$. The
absolute value $|\cdot|_K$ defines a topology on $K$ and turns it
into a locally compact field. We denote by $\mu$ the Haar measure
on $R$, with the usual normalization $\mu(R)=1$. By abuse of
notation, we also denote by $\mu$ the product measure on the
cartesian powers $K^d$ of $K$. If $S$ is a set, we denote its
cardinality by $|S|$. For any complex number $\alpha$, we denote
its real part by $\Re(\alpha)$.

Let $\widetilde{x}=(x_1,\ldots,x_d)$ be a tuple of coordinates,
and let $f$ be a polynomial in $K[\widetilde{x}]$. Let
$\Phi:K^d\rightarrow \C$ be a Schwartz-Bruhat function, i.e. a
locally constant function with compact support, and let
$\chi:R^*\rightarrow \C^*$ be a character of $R^*$, i.e. a group
homomorphism with finite image. We put $\chi(0)=0$. We'll denote
the trivial character by $\chi_{triv}$.

 \begin{definition}[Igusa zeta function]
The Igusa zeta function associated to the triple $(f,\chi,\Phi)$
is defined by
$$Z(f,\chi,\Phi;s)=\int_{K^d}\Phi\cdot(\chi\circ
ac)(f)\cdot|f|_K^s d\mu$$ for $s\in \C$ with $\Re(s)>0$.
\end{definition}

 If $\Phi$ is the characteristic function
$\mathbf{1}_{R^d}$ of $R^d$, or if $\chi$ is the trivial character
$\chi_{triv}$, we'll omit $\Phi$, resp. $\chi$, from the notation.
For instance, we write $Z(f;s)$ instead of
$Z(f,\chi_{triv},\mathbf{1}_{R^d};s)$.
 It is often convenient to introduce a new
variable $t=q^{-s}$; then $Z(f,\chi,\Phi;s)$ is a power series in
$t$ with coefficients in $\C$ which converges for $|t|_K<1$.
Observe that $Z(f,\chi,\Phi;s)$ only depends on $s$ modulo $(2\pi
i/\ln q)\Z$.

The definition of $Z(f,\chi,\Phi;s)$ depends on the choice of the
uniformizer $\pi$ (through the angular component $ac(\cdot)$).
Following Igusa \cite[8.2]{Igusa:intro}, it is more natural to
consider
$$Z(f,\cdot,\Phi;\cdot)$$ as a function on the complex analytic
space $\Omega(K^\ast)$ of quasi-characters $\omega$ of $K^*$
(continuous group homomorphisms $K^*\rightarrow \C^*$). We will
not use this point of view in this article.
%
%
\subsection{$p$-adic cylinders}\label{subsec-pcyl}
Now we introduce the notion of a cylinder in $R^d$. Elementary as
it may be, it will play a central role in the construction of the
motivic measure in Section \ref{subsec-motint}.
\begin{definition}[Cylinder]
 We consider, for each $m\geq 0$, the natural
projection
$$\pi_m:R^d\rightarrow (R/P^{m+1})^d$$
If $S$ is a subset of $(R/P^{m+1})^d$, we call the subset
$(\pi_m)^{-1}(S)$ of $R^d$ the cylinder over $S$.

 We say that a subset $C$ of $R^d$ is a cylinder if it is the
 cylinder over $\pi_m(C)$ for some $m\geq 0$.
\end{definition}
\begin{lemma}\label{lemma-cyl} For any cylinder $C$ in $R^d$, the series
$(q^{-d(m+1)}|\pi_m(C)|)_{m\geq 0}$ is constant for $m\gg 0$, and
its limit is equal to the Haar measure of $C$.

 More precisely, if we choose $m_0\geq 0$ such that $C$ is the cylinder over $\pi_{m_0}(C)$, then for $m\geq
 m_0$ we have
$$q^{-d(m+1)}|\pi_m(C)|=q^{-d(m_0+1)}|\pi_{m_0}(C)|=\mu(C)$$
\end{lemma}
\begin{proof}
For $m\geq m_0$, $C$ can be written as a disjoint union
$$C=\bigsqcup_{a\in \pi_m(C)}(a+(P^{m+1}R)^d)$$ By
translation invariance of the Haar measure and the fact that
$(P^{m+1}R)^d$ has measure $q^{-(m+1)d}$, the measure of $C$
equals
$$q^{-(m+1)d}|\pi_m(C)|$$
\end{proof}
\subsection{Solutions of congruences}\label{subsec-cong}
If $f$ has integer coefficients, i.e. $f\in R[\widetilde{x}]$,
then the local zeta function $Z(f;s)$ is closely related to the
number of solutions of the congruences $f\equiv 0$ modulo powers
of $P$.

\begin{definition}[Igusa Poincar\'e series]
Let $f$ be an element of $R[\widetilde{x}]$. Denote for each $m>
0$ by $N_m(f)$ the cardinality of the set
$$\{z\in (R/P^m)^d\,|\,f(z)= 0 \mbox{ in }R/P^m\}$$
The generating series
$$Q(f;t)=\sum_{m> 0}N_m(f) t^m\in \Z[[t]]$$
is called the Igusa Poincar\'e series associated to $f$.
\end{definition}

 \begin{prop}\label{prop-transform} Putting $t=q^{-s}$, we have
\begin{equation}\label{eq-transf}Q(f;q^{-d}t)=\frac{t}{1-t}(1-Z(f;s))
\end{equation} in $\Q[[t]]$.
\end{prop}
\begin{proof}
 By definition of the zeta
function, the coefficient of $t^i=q^{-is}$ in $Z(f;s)$ is equal to
the Haar measure of the set
$$S_i=\{z\in R^d\,|\,v(f(z))=i\}$$ for each $i\geq 0$. So for each $m> 0$, the
coefficient $c_m$ of $t^m$ in the right hand side of
(\ref{eq-transf}) is $1-\sum_{i=0}^{m-1}\mu(S_i)$. By additivity
of the measure and the fact that $R^d$ has measure $1$, we get
$c_m=\mu(S_{\geq m})$ with $$S_{\geq m}=\{z\in R^d\,|\,v(f(z))\geq
m\})$$ The set $S_{\geq m}$ is a cylinder over the set
$$C_{m-1}:=\{z\in (R/P^m)^d\,|\,f(z)= 0 \mbox{ in }R/P^{m}\}$$
whose cardinality equals $N_m(f)$. By Lemma \ref{lemma-cyl} we
obtain equation (\ref{eq-transf}).\end{proof}
\subsection{Rationality results}\label{subsec-rat}
We fix a polynomial $f\in K[\widetilde{x}]$, a Schwartz-Bruhat
function $\Phi$ on $K^d$
 and a character $\chi$ on $R^*$.
\begin{theorem}[Igusa]\label{thm-rat0}
The local zeta function $Z(f,\chi,\Phi;s)$ is rational in
$t=q^{-s}$.
\end{theorem}
 In particular, $Z(f,\chi,\Phi;s)$ has a meromorphic
continuation to $\C$. The proof actually yields more information:
it specifies a finite subset $S$ of $\Q_{<0}$ such that the poles
of $Z(f,\chi,\Phi;s)$ are contained in $$S+(2\pi i/\ln q)\Z$$
Before giving and proving the precise statement, we need some
additional notation.

The polynomial $f$ defines a morphism of $K$-varieties
$f:\A^d_K\rightarrow \A^1_K$. We choose an \textit{embedded
resolution of singularities} for the morphism $f$. This is a
projective birational morphism of $K$-varieties
$$h:Y\rightarrow \A^d_K$$ such that $Y$ is smooth over $K$, $h$ is an
isomorphism over the complement of the zero locus of $f$ in
$\A^d_K$, and the divisor $E=(f\circ h)$ on $Y$ has strict normal
crossings. Such an embedded resolution always exists, by
\cite{hironaka}.

We denote by $Jac_h$ the Jacobian ideal sheaf of $h$ on $Y$. If
$\omega$ is a gauge form on an open subvariety $U$ of $Y$, then on
$U$ the Jacobian ideal sheaf $Jac_h$ is generated by the unique
regular function $g$ on $U$ satisfying
$$h^*(dx_1\wedge\ldots\wedge dx_d)=g\cdot \omega$$ Since we
assumed that $h$ is an isomorphism on $Y\setminus E$, the divisor
$(Jac_h)$ is supported on $E$.

 We denote by $E_i,\,i\in I$ the irreducible components of
$E$, by $N_i$ the multiplicity of $f\circ h$ along $E_i$, and by
$\nu_i-1$ the multiplicity of the Jacobian ideal $Jac_h$ of $h$
along $E_i$. The couples $(N_i,\nu_i)$ are called the
\textit{numerical data} of the resolution $h$. Now Theorem
\ref{thm-rat0} can be refined as follows.

\begin{theorem}[Igusa]\label{thm-rat1}
The local zeta function $Z(f,\chi,\Phi;s)$ is rational in
$t=q^{-s}$. If $\alpha$ is a pole of $Z(f,\chi,\Phi;s)$, then
there exists an index $i\in I$ such that the order of $\chi$
divides $N_i$ and such that $\alpha$ is contained in
$$-\frac{\nu_i}{N_i}+\frac{2\pi i}{\ln q}\Z$$
\end{theorem}
\begin{proof}
For the proof of the theorem we assume that the reader is familiar
with the theory of integration of differential forms $\omega$ of
maximal degree on $K$-analytic manifolds $M$
\cite[7.4]{Igusa:intro}. Recall in particular that $\omega$
defines a measure on the set of compact open subsets of $M$. Since
$Y$ is a smooth $K$-variety, we can view its set of rational
points $Y(K)$ as a $K$-analytic manifold. The change of variables
formula for $p$-adic integrals \cite[7.4.1]{Igusa:intro} yields
$$Z(f,\chi,\Phi;s)=\int_{Y(K)}(\Phi\circ h)\cdot(\chi\circ ac)(f\circ h)\cdot
|f\circ h|_K^s h^*(d\widetilde{x})$$ with
$d\widetilde{x}=dx_1\wedge\ldots\wedge dx_d$, and this integral
can be computed locally on $Y(K)$ as follows.

Let $b$ be a point of $Y(K)$ and let $\{i_1,\ldots,i_r\}$ be the
(possibly empty) set of indices $i$ in $I$ with $b\in E_i$. Since
the divisor $E$ has strict normal crossings, there exist an open
neighbourhood $U$ of $b$ in $Y(K)$, and analytic coordinates
$y_1,\ldots,y_d$ and nowhere-vanishing $K$-analytic functions $u$,
$v$ on $U$, such that
$$f\circ h=u\prod_{j=1}^r (y_j)^{N_{i_j}}\mbox{ and
}h^*d\widetilde{x}=v\prod_{j=1}^r(y_j)^{\nu_{i_j}-1}dy$$ Since the
value group $|K^*|_K$ is discrete, we may even assume that $|u|_K$
and $|v|_K$ are constant on $U$.

Since $h$ is proper, the support of $\Phi \circ h$ is compact, so
we can write $Z(f,\chi,\Phi;s)$ as a finite $\C$-linear
combination of finite products with factors of the form
$q^{as}\cdot A(e,i)$. Here $e,a\in \Z$, $i\in I$, and
$$A(e,i):=\int_{z\in \pi^{e}R}\chi^{N_i}(ac(z))\cdot |z|_K^{N_is+\nu_i-1}dz$$

 If $\chi^{N_i}\neq 1$ then $A(e,i)$ vanishes, by a standard
argument: if $w\in R^*$ then a substitution $z'=wz$ shows that
$$A(e,i)=\chi^{N_i}(w)\cdot A(e,i)$$
If $\chi^{N_i}= 1$
then
\begin{eqnarray*}
A(e,i)&=&\sum_{m\geq e}\mu(\{z\in R\,|\,v(z)=m\})\cdot
q^{-m(N_is+\nu_i-1)}
\\&=&(1-\frac{1}{q})\frac{q^{-e(N_is+\nu_i)}}{1-q^{-(N_is+\nu_i)}}
\end{eqnarray*}
\end{proof}
\begin{cor} If $f$ has coefficients in $R$, the Igusa Poincar\'e series
$Q(f;T)$ is rational over $\Q$.
\end{cor}
\begin{proof}
 This follows from Proposition \ref{prop-transform}. \end{proof}
 The rationality of $Q(f;t)$
 was conjectured by Borevich and Shafarevich. Denef gave a
vast generalization of Igusa's rationality result in Theorem
\ref{thm-rat0}, avoiding resolution of singularities but instead
using $p$-adic cell decomposition \cite{Denef} to prove a general
structure theorem on definable $p$-adic integrals. As a special
case, Denef obtained the rationality of the series
$$\widetilde{Q}(f;t)=\sum_{m>0}\widetilde{N}_m(f)t^m$$ where $\widetilde{N}_m(f)$
denotes the cardinality of the image of the projection
$$\{x\in R^d\,|\,f(x)=0\} \rightarrow (R/P^m)^d$$
 The
rationality of this series was conjectured by Serre and
Oesterl\'e. To my knowledge, there is still no ``purely''
arithmetic or geometric proof which does not use $p$-adic model
theory. If the reduction $\overline{f}$ of $f$ modulo $P$ is
non-zero and defines a smooth hypersurface in $\A^d_k$, then
$\widetilde{Q}(f;t)=Q(f;t)$, but in the singular case the
geometric meaning of its poles is completely unknown. See
\cite{DL}\cite{LejReg}\cite{Nic2}\cite{Nic1} for some particular
cases in the motivic setting.

Igusa's theorem yields a complete set of candidate poles for
$Z(f,\chi,\Phi;s)$, but many of these will not be actual poles of
the local zeta function. For one thing, the set of candidate poles
depends on the embedded resolution $h$. But even if we take the
intersection of these sets over all resolutions $h$, the resulting
set will in general still be too big.

 For instance, in the curve case $d=2$, there exists a
\textit{minimal} embedded resolution, but not all the
corresponding candidate poles will be poles of the zeta function.
This phenomenon is related to the monodromy conjecture, which puts
additional (conjectural, if $d>2$) restrictions on the poles.

\subsection{Denef's formula}
If the embedded resolution $h$ has good reduction modulo the
maximal ideal $P$ of $R$ (in a certain technical sense; see below)
then Denef gave
a very explicit formula for $Z(f,\chi,\Phi;s)$ in terms of the
resolution $h$ \cite{Denef5}\cite{Denef6}.
 If $f$ and $h$ are defined over a number field $L$, then
good reduction holds at almost all finite places $\mathfrak{P}$ of
$L$.

Denef's formula involves the numerical data of $h$ and certain
character sums over sets of rational points on the reductions of
the exceptional components $E_i$ modulo $P$. Using the \'etale
Kummer sheaf associated to $\chi$ and Grothendieck's trace
formula, this yields a cohomological interpretation for the local
zeta function \cite{Denef6}.

To state Denef's formula, we need some additional notation. A
Schwartz-Bruhat function $\Phi$ on $K^d$ is called residual if its
support is contained in $R^d$ and the value of $\Phi$ at a point
$y$ of $R^d$ only depends on its residue class in $k^d$. Let $f$
be an element of $R[\widetilde{x}]$, and $h:Y\rightarrow \A^d_K$
an embedded resolution of singularities for $f$.
%
 We fix
a closed immersion of $Y$ into projective space $\Pro^r_{\A^d_K}$
for some $r>0$. We denote by $E=(f\circ h)$ the divisor on $Y$
defined by $f\circ h$, by $E_i,\,i\in I$ its irreducible
components, and by $(N_i,\nu_i)$ the corresponding numerical data.
For any integer $e>0$, we denote by $\mathcal{I}_e$ the set of
indices $i\in I$ such that $e|N_i$.

We denote by $\overline{f}$ the reduction of $f$ modulo $P$. It is
an element of $k[\widetilde{x}]$. For any closed subvariety $Z$ of
$Y$ we denote by $\overline{Z}$ the reduction modulo $P$ of the
Zariski-closure of $Z$ in $\Pro^r_{\A^d_{R}}$. The morphism $h$
induces a morphism of $k$-varieties
$\overline{h}:\overline{Y}\rightarrow \A^d_{k}$.

For any subset $J$ of $I$ we put $E_J=\cap_{j\in J}E_j$ and
$\overline{E}^o_J=\overline{E_J}\setminus \cup_{i\notin
J}\overline{E}_i$. In particular, $E_{\emptyset}=Y$ and
$\overline{E}^o_{\emptyset}=\overline{Y}\setminus \overline{E}$.
We put $m_J=gcd\{N_i\,|\,i\in J\}$. Note that
$m_{\emptyset}=gcd(\emptyset)=0$.
 For $i\in I$
we write  $\overline{E}^o_i$ instead of
 $\overline{E}^o_{\{i\}}$.

Following Denef, we say that the resolution $h$ has good reduction
modulo $P$ if $\overline{Y}$ and all $\overline{E}_i$ are smooth,
$\overline{E}$ is a divisor with strict normal crossings,  and
$\overline{E}_i$ and $\overline{E}_j$ have no common components if
$i\neq j$. In this case, $\overline{E_J}=\cap_{j\in
J}\overline{E}_j$ for each subset $J$ of $I$. We say that the
resolution has tame good reduction if, moreover, none of the
multiplicities $N_i$ belongs to $P$.  If $h$ has good reduction,
$J$ is a subset of $I$, and $a$ is a point of
$\overline{E}_J^o(k)$, then in the local ring
$\mathcal{O}_{\overline{Y},a}$ the element $\overline{f}\circ
\overline{h}$ can be written as $u\cdot v^{m_J}$ with $u$ a unit.
If $\chi$ is a character of $R^*$ which is trivial on $1+P$, then
$\chi(u(a))$ is well-defined since $u(a)\in k^*=R^*/(1+P)$. If,
moreover, $e$ is an element of $\Z_{>0}$ such that $J\subset
\mathcal{I}_e$ and $\chi^e$ is trivial, then $e|m_J$ and
$$\omega_{\chi}(a):=\chi(u(a))$$ does not depend on the choice of
$u$ and $v$.

 If $\Phi$ is a residual Schwartz-Bruhat function on $K^d$, we denote by
$\overline{\Phi}$ the induced function $k^d\rightarrow \C$. Let
$\chi$ be a character of $R^*$, of order $e$. If $\chi$ is trivial
on $1+P$ and $h$ has good reduction modulo $P$, then for any
subset $J$ of $\mathcal{I}_e$ we put
$$
c_{J,\Phi,\chi}=\sum_{a\in\overline{E}_J^o(k)}(\overline{\Phi}\circ
\overline{h})(a)\cdot \omega_{\chi}(a)$$ In particular, if $\Phi$
is the characteristic function of $R^d$ and $\chi$ is trivial,
then
$$c_{J,\Phi,\chi}=|\overline{E}_J^o(k)|$$ for all
subsets $J$ of $I$.

\begin{theorem}[Denef]
\label{theo-denef} Assume that $\overline{f}\neq 0$. Let $\Phi$ be
a residual Schwartz-Bruhat function on $K^d$ and $\chi$ a
character of $R^*$ of order $e$. If $\chi$ is non-trivial on $1+P$
and $h$ has tame good reduction, then $Z(f,\chi,\Phi;s)$ is
constant as a function of $s$. If $\chi$ is trivial on $1+P$ and
$h$ has good reduction, then we have
$$Z(f,\chi,\Phi;s)=q^{-d}\sum_{J\subset
\mathcal{I}_e}\left(c_{J,\Phi,\chi}\prod_{j\in
J}\frac{(q-1)q^{-N_js-\nu_j}}{1-q^{-N_js-\nu_j}}\right)$$
\end{theorem}

If $f$ and $h$ are defined over a number field $L$, then $h$ has
tame good reduction modulo $\mathfrak{P}$ for almost all finite
places $\mathfrak{P}$ of $L$ (i.e. all but a finite number).
Hence, if we denote by $L_{\mathfrak{P}}$ the completion of $L$ at
$\mathfrak{P}$ and by $\mathcal{O}_{\mathfrak{P}}$ its ring of
integers, then for all but a finite number of $\mathfrak{P}$,
Theorem \ref{theo-denef} gives an explicit expression for the
$\mathfrak{P}$-adic zeta function associated to $f$, any residual
Schwartz-Bruhat function $\Phi$ on $(L_{\mathfrak{P}})^d$ and any
character $\chi$ of $\mathcal{O}_{\mathfrak{P}}^*$.
\section{The $p$-adic monodromy
conjecture}\label{sec-pmoncon}
\subsection{The Milnor fibration}\label{subsec-milnor} Let $X$ be a
complex manifold, and let
 $$g:X\rightarrow \C$$ be a non-constant analytic map. We denote by $X_s$ the special fiber of $g$ (i.e.
 the analytic space defined by $g=0$), and we fix a
  point $x\in X_s$.

What does
 $X_s$ look like in a neighbourhood of $x$?
 If $g$ is smooth at $x$ this is easy: $X_s$ is locally a
 complex submanifold of $X$. If $g$ is not smooth at $x$, then the topology of $X_s$ near $x$ can be
 studied by means of the Milnor fibration \cite{Milnor}\cite{Dimca}.

Working locally, we may assume that $X=\C^{d}$. Let
$B=B(x,\varepsilon)$ be an open ball around $x$ in $\C^{d}$ with
radius $\varepsilon$, let $D=D(0,\eta)$ be an open disc around the
origin $0$ in $\C$ with radius $\eta$, and put
$D^\ast=D\setminus\{0\}$. For $0<\eta \ll \varepsilon\ll 1$ the
map
$$g_x:g^{-1}(D^*)\cap B\rightarrow D^\ast$$ is a locally
trivial fibration, called the Milnor fibration of $g$ at $x$.

We consider the universal covering space
$$\widetilde{D^*}=\{z\in \C\,|\,\Im(z)>-\log \eta\}\rightarrow D^*:z\mapsto \exp(iz)$$
of $D^*$ and we put $$F_x=(g^{-1}(D^*)\cap
B)\times_{D^\ast}\widetilde{D^*}$$ This is the universal fiber of
the fibration $g_x$, and it is called the Milnor fiber of $g$ at
$x$. Since $g_x$ is a locally trivial fibration and
$\widetilde{D^*}$ is contractible, $F_x$ is homotopy-equivalent to
the fiber of $g_x$ over any point of $D^*$.

If $g$ is smooth at $x$, then the fibration $g_x$ is trivial. In
general, the defect of triviality is measured by the monodromy
action on the singular cohomology of $F_x$, i.e. the action of the
group $\pi_1(D^*)$ of covering transformations of
$\widetilde{D^*}$ over $D^*$ on $\oplus_{i\geq
0}H^i_{sing}(F_x,\Z)$.
 The action of the canonical generator $z\mapsto z+2\pi$ is
called the monodromy transformation and denoted by $M_x$. We say
that a complex number $\gamma$ is a monodromy eigenvalue of $g$ at
$x$ if $\gamma$ is an eigenvalue of the monodromy transformation
$M_x$ on $H_{sing}^i(F_x,\Z)$ for some $i\geq 0$. These monodromy
eigenvalues are roots of unity.
\subsection{The monodromy conjecture}\label{subsec-pmoncon}
 Now assume that $f$
is a polynomial in the variables $\widetilde{x}=(x_1,\ldots,x_d)$
over some number field $L$. Then, for any finite place
$\mathfrak{P}$, we can view $f$ as a polynomial over the
$\mathfrak{P}$-adic completion $L_{\mathfrak{P}}$ of $L$ and
consider the associated local zeta functions
$Z_\mathfrak{P}(f,\chi,\Phi;s)$ for varying $\Phi$ and $\chi$.  On
the other hand, we can view $f$ as a complex polynomial, defining
an analytic map
$$f:X=\C^d\rightarrow \C$$ and we can consider the
singularities of its special fiber $X_s$. These objects are
related by Igusa's Monodromy Conjecture. We denote by
$\mathcal{O}_{\mathfrak{P}}$ the ring of integers of
$L_{\mathfrak{P}}$.

\begin{conjecture}[Monodromy Conjecture]\label{con-mon} Let $L$ be a number field and $f$ an element
of $L[\widetilde{x}]$. For almost all finite places $\mathfrak{P}$
of $L$, we have the following property:

if $\Phi$ is a Schwartz-Bruhat function on $L_{\mathfrak{P}}$,
$\chi$ is a character of $\mathcal{O}_{\mathfrak{P}}^*$, and
$\alpha$ is a pole of the local zeta function
$Z_\mathfrak{P}(f,\chi,\Phi;s)$, then $\exp(2\pi i\Re(\alpha))$ is
a monodromy eigenvalue of $f$ at some point $x\in X_s$.
\end{conjecture}
 There is a stronger version of the conjecture, saying that
under the same conditions, $\Re(\alpha)$ is a root of the
Bernstein-Sato polynomial $b_f(s)$ of $f$ ; see \cite[\S
4]{Igusa:intro} for the definition of $b_f(s)$. This statement is
indeed stronger since it is known that for any root $\beta$ of
$b_f(s)$, the value $\exp(2\pi i \beta)$ is an eigenvalue of
monodromy at some point $x$ of $X_s$ \cite{Mal2}\cite{kashi2}. It
is also known that the roots of $b_f(s)$ are negative rational
numbers \cite{kashi1}. If $X_s$ is smooth, then $b_f(s)=s+1$, so
the strong version of the monodromy conjecture is valid in this
(very) special case.

For future reference, we state the following particular case of
the Monodromy Conjecture.
\begin{conjecture}[Untwisted Monodromy Conjecture]\label{con-untwist}
Let $L$ be a number field and $f$ an element of
$L[\widetilde{x}]$. For almost all finite places $\mathfrak{P}$,
we have the following property:

if  $\alpha$ is a pole of the local zeta function
$Z_\mathfrak{P}(f;s)$ associated to $f$ over $L_{\mathfrak{P}}$
(and to $\Phi=1_{\mathcal{O}_{\mathfrak{P}}^d}$ and
$\chi=\chi_{triv}$), then $\exp(2\pi i\Re(\alpha))$ is a monodromy
eigenvalue of $f$ at some point $x\in X_s$.
\end{conjecture}
\subsection{Some evidence}\label{subsec-evidence}
(1) \textit{The archimedean case.} As we mentioned in the
introduction, Igusa's theory of local zeta functions also has a
archimedean wing, studying zeta functions over the local field
$\mathbb{K}$ with $\mathbb{K}=\R$ or $\mathbb{K}=\C$. They are
defined in a similar way, as functions on the space of
quasi-characters of $\mathbb{K}$, associated to a  polynomial $f$
over $\mathbb{K}$ and a Schwartz-Bruhat function $\Phi$ on
$\mathbb{K}^d$ (a $\mathcal{C}^{\infty}$-function with compact
support).

 For instance, for $\mathbb{K}=\R$ and $f\in \R[\widetilde{x}]$ with $\widetilde{x}=(x_1,\ldots,x_d)$, we get
$$Z(f,\Phi,\chi;s)=\int_{\R^d}\Phi\cdot \chi(f)\cdot |f|^sd\widetilde{x}$$
for $s\in \C$ with $\Re(s)>0$, where $\chi$ is either the constant
function $1$ or the sign function $sgn(\cdot)$.

 Using the functional equation
$D\cdot f^{s+1}=b_f(s)f^s$ for the Bernstein-Sato polynomial (with
$D\in K[\widetilde{x},\partial_{\widetilde{x}},s]$) and
integration by parts, it is not hard to show that $Z(f,\Phi;s)$
has a meromorphic continuation to $\C$, and that its poles
$\alpha$ are of the form $\beta-j$ with $b_f(\beta)=0$ and $j\in
\N$ \cite[5.3]{Igusa:intro}. In particular, $\exp(2\pi i\alpha)$
is an eigenvalue of monodromy\footnote{In spite of this result, it
remains a challenging problem to determine which candidate poles
actually occur, and to determine their multiplicity; see for
instance \cite{DeNiSa}.}.

 However, integration by parts does not make sense
in the $p$-adic setting. Therefore, it is quite surprising that we
still get (at least conjecturally) a similar relation between
poles of the zeta function and roots of the Bernstein polynomial.

(2) \textit{A'Campo's formula.} Let $g$ be a polynomial in
$\C[x_1,\ldots,x_d]$ and denote by $Y_s$ the complex hypersurface
defined by $g$. The eigenvalues of monodromy of $g$ at the points
of $Y_s(\C)$ can be computed on an embedded resolution of
singularities $h:Y\rightarrow \A^d_{\C}$ of $g$. Denote by
$E_i,\,i\in I$ the irreducible components of the divisor
$E=(g\circ h)$ and by $(N_i,\nu_i)$ the corresponding numerical
data. For each $i\in I$ we put $E_i^o=E_i\setminus \cup_{j\neq
i}E_j$.

If $x$ is a point of $Y_s(\C)$, then the monodromy zeta function
$\zeta_{g,x}(T)$ is defined as the alternating product of the
characteristic polynomials of the monodromy transformation on the
singular cohomology spaces of the Milnor fiber $F_x$ of $g$ at
$x$~:
$$\zeta_{g,x}(T)=\prod_{i\geq
0}\mathrm{det}(1-T\cdot M_x\,|\,H^i_{sing}(F_x,\Z))^{(-1)^{i+1}}$$
 Using Leray's spectral sequence and an explicit
description of the stalks of the complex of nearby cycles of
$g\circ h$, A'Campo \cite{A'C} proved the following formula:
$$\zeta_{g,x}(T)=\prod_{i\in I}(1-T^{N_i})^{-\chi_{top}(E_i^o\cap
h^{-1}(x))}$$ Here $\chi_{top}(.)$ is the topological Euler
characteristic.

 Moreover, using the perversity of the nearby cycles
complex, Denef \cite{Denef5} observed that any eigenvalue of
monodromy at a point $x\in Y_s(\C)$ occurs as a zero or pole of
$\zeta_{g,y}(T)$ for some (possibly different) point $y\in
Y_s(\C)$.

Now let $L$ be a number field, and $f$ a polynomial in
$L[x_1,\ldots,x_d]$, and fix an embedded resolution of
singularities for $f$. In principle, using Denef's explicit
formula (Theorem \ref{theo-denef}) and A'Campo's formula, one can
compute on this embedded resolution the residues of the candidate
poles of the $\mathfrak{P}$-adic zeta function of $f$ (for almost
all finite places $\mathfrak{P}$ of $L$), and all eigenvalues of
monodromy. By this procedure, one eliminates fake candidate poles
and one obtains a list of eigenvalues of monodromy, and when all
remaining candidate poles induce eigenvalues, the monodromy
conjecture is proven. For instance, using A'Campo's formula, one
easily sees that $1$ and $\exp(\pi i/3)$ are monodromy eigenvalues
of the polynomial $g(x,y)=y^2-x^3$ at the origin. Hence, the
monodromy conjecture holds for $g(x,y)$, since the computations in
Section \ref{sec-intro} show that for $p\gg 0$ the poles of
$Z_p(g,s)$ have real part $-1$ or $-5/6$.

 Unfortunately, this strategy requires a strong
control over the configuration of the exceptional components $E_i$
and their numerical data (one has to prove that configurations
inducing ``bad'' poles cannot occur in the exceptional locus of a
resolution).
 It led, for instance, to a proof of the conjecture in the
 following cases:
\begin{itemize}
\item $d=2$ \cite{Loepadic}\cite{Rod} \item $d=3$ and $f$
homogeneous
 \cite{rove-holo} \cite{aclm1} \item
superisolated surface singularities
 \cite{aclm1}
\end{itemize} More generally, Veys obtained very nice results on
possible configurations of exceptional divisors, thus gathering
strong evidence for the conjecture, especially for $d=3$
\cite{Veysnum}\cite{Veys-crelle}.

 In other settings, some
combinatorial description of the singularities yields an
expression for the local zeta function and/or the monodromy zeta
function, and the conjecture can be proven from this expression;
e.g. for
\begin{itemize}
\item $f$ non-degenerate w.r.t. its Newton polyhedron, with an
additional technical assumption \cite{loe-newton} \item $f$
quasi-ordinary \cite{aclm2}
\end{itemize}

(3) \textit{Prehomogeneous vector spaces.} Finally, we should
mention the case where $f$ is a relative invariant of a
prehomogeneous vector space $(G,X)$. This was one of the first
classes of examples where Igusa made extensive computations
(exploiting the group structure), and these computations led to
several conjectures, including the monodromy conjecture.
 The monodromy conjecture was proven in \cite{kisazhu}
for $(G,X)$ irreducible and reduced, using Igusa's group-theoretic
expression for a complete list of candidate-poles of the local
zeta function.

\bigskip
 It does not seem probable that any of these techniques can
be used to deal with the general case, because the geometric
complexity becomes quite hard to control in higher dimensions; a
more intrinsic relation between the local zeta function and
monodromy should be discovered.

 Recently, Sebag and the author introduced a geometric object whose
rational points are related to the (motivic) zeta function, and
whose \'etale cohomology (with Galois action) coincides with the
cohomology of the Milnor fiber (with monodromy action)
\cite{NiSe}\cite{Ni-trace}. This object is a non-archimedean
analytic space over the field $\C((t))$ of complex Laurent series;
we called it the analytic Milnor fiber. Its geometric/arithmetic
properties are closely related to the structure of the
singularities of $f$, and many of its invariants have a natural
interpretation in singularity theory. We hope that the study of
this object will lead to new insights into the monodromy
conjecture.

\section{Motivic integration}\label{sec-motint}
\subsection{From $\Z_p$ to $k{[[}t{]]}$}\label{subsec-motintro} 

\noindent Kontsevich invented motivic integration to strengthen
the following result by Batyrev \cite{Baty}.
\begin{theorem}[Batyrev] If two complex Calabi-Yau varieties are
birationally equivalent, then they have the same Betti
numbers.\end{theorem}
 Batyrev proved this result using $p$-adic integration and
the Weil Conjectures. Kontsevich observed that Batyrev's proof
could be ``geometrized'', avoiding the passage to finite fields
and yielding a stronger result: equality of Hodge numbers. The key
was to replace the $p$-adic integers by $\C[[t]]$, and $p$-adic
integration by motivic integration.

Kontsevich presented these ideas at a famous ``Lecture at Orsay''
in 1995, but never published them. The theory was developed and
generalized in the following directions:
\begin{itemize}
 \item Denef and Loeser \cite{DLinvent} developed a theory of \textit{geometric}
motivic integration on arbitrary algebraic varieties over a field
of characteristic zero.
 They also created a theory of \textit{arithmetic} motivic
integration \cite{DL}, with good specialization properties to
$p$-adic integrals in a general setting, using the model theory of
pseudo-finite fields. The motivic integral appears here as a
universal integral, specializing to the $p$-adic ones for almost
all $p$.

 \item Loeser and Sebag constructed a theory of motivic
integration on formal schemes \cite{sebag1} and rigid varieties
\cite{motrigid}, working over an arbitrary complete discrete
valuation ring with perfect residue field.

\item Cluckers and Loeser built a very general framework for
motivic integration theories, based on model theory \cite{cluloe}.
A different model-theoretic approach was developed by Hrushovksi
and Kazhdan \cite{hruka}.
\end{itemize}
 We will only discuss the so-called ``na\"ive'' geometric
motivic integration on smooth algebraic varieties (but the
adjective ``na\" ive'' does no right to the stunning vision behind
the constructions). We start by explaining the basic ideas.

In Section \ref{subsec-def}, we introduced the notion of cylinder
in $(\Z_p)^d$ using the projection maps $\pi_m:(\Z_p)^d\rightarrow
(\Z/p^{m+1})^d$ for $m\geq 0$, and we saw that the Haar measure of
a cylinder $C$ can be computed from the cardinality of the
projection $\pi_m(C)$ for $m\gg 0$.
 If we identify $\Z_p$ with the ring of Witt vectors
$W(\mathbb{F}_p)$, then the map $\pi_m$ simply corresponds to the
truncation map
$$W(\mathbb{F}_p)\rightarrow W_{m+1}(\mathbb{F}_p):(a_0,a_1,\ldots)\mapsto
(a_0,a_1,\ldots,a_m)$$

 The idea behind the theory of motivic integration is to
make a similar construction, replacing $W(\mathbb{F}_p)$ by the
ring of formal power series $k[[t]]$ over some field $k$, and the
map $\pi_m$ by the truncation map
$$k[[t]]^d\rightarrow (k[t]/(t^{m+1}))^d:(\sum_{i\geq 0}a_{1,i}t^i,\ldots,\sum_{i\geq 0}a_{n,i}t^i)\mapsto
(\sum_{i=0}^{m}a_{1,i} t^i,\ldots,\sum_{i=0}^{m}a_{n,i} t^i)$$

 The problem is to give
meaning to the expression $|\pi_m(C)|$ if $C$ is a ``cylinder'' in
$k[[t]]^d$ for infinite fields $k$, and to find a candidate to
replace $p$.
 But interpreting the coefficients of a power series
as affine coordinates, the set $(k[[t]]/(t^{m+1}))^d$ gets the
structure of the set of $k$-points on an affine space
$\A^{(m+1)d}_k$, and if we restrict to cylinders $C$ such that
$\pi_m(C)$ is constructible in $\A^{(m+1)d}_k$, we can use the
Grothendieck ring of varieties as a universal way to ``count''
points on constructible subsets of an algebraic variety.
 The cardinality $p$ of $\mathbb{F}_p$ is
replaced by the ``number'' of points on the affine line $\A^1_k$;
this is the \textit{Lefschetz motive}  $\mathbb{L}$.

The price to pay is that we leave classical integration theory
since our value ring will be an abstract object (the Grothendieck
ring of varieties) instead of $\R$.

In the following sections, we will make these ideas precise.

\subsection{The Grothendieck ring of
varieties}\label{subsec-groth}
Let $k$ be any field. A $k$-variety is a reduced separated
$k$-scheme of finite type.

\begin{definition}[Grothendieck ring of $k$-varieties]
As an abelian group, $K_0(Var_k)$ has the following presentation:
\begin{itemize}
\item generators: isomorphism classes $[X]$ of separated
$k$-schemes of finite type $X$
 \item relations: $[X]=[X\setminus Y]+[Y]$ if $Y$ is a closed
 subscheme of
$X$ (``scissor relations'')
\end{itemize}
We endow $K_0(Var_k)$ with the unique ring multiplication  such
that $[X_1]\cdot[X_2]=[X_1\times_k X_2]$ for any pair $X_1,\,X_2$
of separated $k$-schemes of finite type.
\end{definition}

We put $\LL=[\A^1_k]$ and we denote by $\mathcal{M}_k$ the
localized Grothendieck ring $$\mathcal{M}_k=K_0(Var_k)[\LL^{-1}]$$
The Grothendieck ring and its localization are still quite
mysterious. It is known that $K_0(Var_k)$ is not a domain if $k$
is a field of characteristic zero \cite{Poo}. It is not known if
the localization morphism $K_0(Var_k)\rightarrow \mathcal{M}_k$ is
injective, i.e. if $\LL$ is a zero divisor in $K_0(Var_k)$, or if
$\mathcal{M}_{k}$ is a domain if $k$ is algebraically closed. For
related questions and results, see for instance
\cite{lalu-julienliu} and \cite{Ni-tracevar}.

\begin{remark}
If $X$ is any $k$-scheme of finite type and $X_{red}$ is the
maximal reduced closed subscheme of $X$, then the closed immersion
$X_{red}\rightarrow X$ is a bijection and the scissor relations
imply that $[X]=[X_{red}]$ in $K_0(Var_k)$. Hence, we get the same
Grothendieck ring if we replace ``separated $k$-scheme of finite
type'' by ``$k$-variety'', and the reader who is unfamiliar with
the formalism of schemes can stick to this definition. One word of
warning, though: if $k$ is imperfect and $X_1$ and $X_2$ are
$k$-varieties, $X_1\times_k X_2$ need not be reduced.
\end{remark}

A subset of a  $k$-variety $X$ is called \textit{locally closed}
if it is open in its closure in $X$ w.r.t. the Zariski topology on
$X$. Such a subset carries a unique structure of subvariety of
$X$. A subset $C$ of $X$ is called \textit{constructible} if it is
a finite union of locally closed subsets of $X$. Then we can
always write $C$ as a finite disjoint union of locally closed
subsets $U_1,\ldots,U_r$ and the scissor relations in the
Grothendieck ring imply that the class
$$[C]:=\sum_{i=1}^{r}[U_i]$$ in $K_0(Var_k)$ does not depend on
the choice of the partition.

Let $F$ be a $k$-variety. A morphism of $k$-varieties
$Y\rightarrow X$ is a Zariski-locally trivial fibration with fiber
$F$ if any point of $X$ has an open neighbourhood $U$ such that
$Y\times_X U$ is isomorphic to $F\times_k U$ as a $U$-scheme. In
this case, we have $[Y]=[X]\cdot[F]$ in $K_0(Var_k)$. Indeed,
using the scissor relations and Noetherian induction we can reduce
to the case where the fibration is trivial.


By its very definition, the Grothendieck ring is the universal
additive multiplicative invariant of $k$-varieties: whenever
$\chi(\cdot)$ is an invariant of $k$-varieties taking values in a
ring $A$ and such that $\chi(\cdot)$ is additive w.r.t. closed
immersions and multiplicative w.r.t. the fiber product over $k$,
it will factor through a unique morphism of rings
$$\chi:K_0(Var_k)\rightarrow A$$ with $\chi([X])=\chi(X)$ for any
$k$-variety $X$. So in a way, taking the class $[X]$ of $X$ in the
Grothendieck ring is the most general way to ``count points'' on,
or ``measure the size'' of, the variety $X$.

Here are some important specialization morphisms.

\vspace{5pt} (1) For any finite field $\mathbb{F}_q$, consider the
invariant $\sharp$ which associates to a $\mathbb{F}_q$-variety
$X$ the number of $\mathbb{F}_q$-rational points on $X$. This
invariant is additive and multiplicative and hence defines a ring
morphism $\sharp:K_0(Var_{\mathbb{F}_q})\rightarrow \Z$. It
localizes to a ring morphism
$\sharp:\mathcal{M}_{\mathbb{F}_q}\rightarrow \Z[q^{-1}]$.

\vspace{5pt} (2) For $k=\C$, we can consider the invariant
$\chi_{top}$ which associates to a $\C$-variety $X$ the
topological Euler characteristic of $X(\C)$ w.r.t. the complex
topology. Again, this invariant defines a ring morphism
$\chi_{top}:K_0(Var_{\C})\rightarrow \Z $ which localizes to a
ring morphism $\chi_{top}:\mathcal{M}_{\C}\rightarrow \Z $ since
$\chi_{top}(\A^1_{\C})=1$.

 If $k$ is any field we can
consider the $\ell$-adic Euler characteristic $\chi_{top}$
instead, with $\ell$ a prime invertible in $k$. It is known that
$\chi_{top}$ does not depend on the choice of $\ell$.

\vspace{5pt} (3) For $k=\C$, we can consider the Hodge-Deligne
polynomial $HD(X;u,v)\in \Z[u,v]$ of a $\C$-variety $X$. It is
defined by
$$HD(X;u,v)=\sum_{p,q\geq 0}\sum_{i\geq 0}(-1)^i h^{p,q}(H^i_c(X(\C),\C))u^pv^q$$
where $h^{p,q}(H^i_c(X(\C),\C))$ denotes the dimension of the
$(p,q)$-component of Deligne's mixed Hodge structure on
$H^i_c(X(\C),\C)$. One can show that $HD(\cdot;u,v)$ is additive
and multiplicative, so there exists a unique ring morphism
$$HD:K_0(Var_{\C})\rightarrow \Z[u,v]$$ mapping $[X]$ to
$HD(X;u,v)$ for each complex variety $X$. We have
$$HD(\A^1_{\C};u,v)=uv$$ and $HD$ localizes to a ring morphism
$$HD:K_0(Var_{\C})\rightarrow \Z[u,u^{-1},v,v^{-1}]$$

\subsection{Arc spaces}\label{subsec-arc}
 Let $X$ be a variety over $k$. For each integer $n\geq
0$,
 we consider the functor
 $$F_n:(k-alg)\rightarrow (Sets):A\mapsto X(A[t]/(t^{n+1}))$$
 which sends a $k$-algebra $A$ to the set of points on $X$ with
 coordinates in $A[t]/(t^{n+1})$.

 \begin{prop}\label{prop-arc} The functor $F_n$ is representable by a separated $k$-scheme of finite
 type $\mathcal{L}_n(X)$, called the $n$-th jet scheme of $X$. If
 $X$ is affine, then so is $\mathcal{L}_n(X)$.\end{prop}

 The proposition means that there exists for any $k$-algebra $A$
 a bijection $\phi_n(A)$ between the set of points on $\mathcal{L}_n(X)$ with
 coordinates in $A$ and the set $F_n(A)$, with the property that
 the diagram
$$\begin{CD}
\mathcal{L}_n(X)(A)@>\phi_n(A)>> F_n(A)
\\@VVV @VVV
\\\mathcal{L}_n(X)(B)@>\phi_n(B)>> F_n(B)
\end{CD}$$
 commutes for any morphism of $k$-algebras $A\rightarrow B$. By Yoneda's Lemma, this
 property determines $\mathcal{L}_n(X)$ as a $k$-scheme, up to canonical isomorphism.

Instead of giving a proof of Proposition \ref{prop-arc}, we look
at an example.

\begin{example}\label{ex}
Let $X$ be the closed subvariety of $\A^2_k=\Spec k[x,y]$ defined
by the equation $x^2-y^3=0$. Then a point of $\mathcal{L}_2(X)$
with coordinates in some $k$-algebra $A$ is a couple
$$(x_0+x_1t+x_2t^2,y_0+y_1t+y_2t^2)$$ with $x_0,\ldots,y_2\in A$
such that
$$(x_0+x_1t+x_2t^2)^2-(y_0+y_1t+y_2t^2)^3\equiv 0\ \mod t^3 $$

Working this out, we get the equations $$\left\{
\begin{array}{lll}
(x_0)^2-(y_0)^3&=&0
\\ 2x_0 x_1-3(y_0)^2y_1&=&0
\\ (x_1)^2+2x_0x_2-3y_0(y_1)^2-3(y_0)^2y_2&=&0
\end{array}\right.$$
and if we view $x_0,\ldots,y_2$ as affine coordinates, these
equations define $\mathcal{L}_2(X)$ as a closed subscheme of
$\A^6_k$.
\end{example}

For any pair of integers $m\geq n\geq 0$ and any $k$-algebra $A$,
the truncation map
$$A[t]/t^{m+1}\rightarrow A[t]/t^{n+1}$$ defines a natural
transformation of functors $F_m\rightarrow F_n$, so by Yoneda's
Lemma we get a natural truncation morphism of $k$-schemes
$$\pi^m_n:\mathcal{L}_m(X)\rightarrow \mathcal{L}_n(X)$$
This is the unique morphism such that for any $k$-algebra $A$, the
square
$$\begin{CD}
X(A[t]/t^{m+1}) @>>> X(A[t]/t^{n+1})
\\ @V\phi_m(A) VV @VV\phi_n(A) V
\\ \mathcal{L}_m(X)(A)@>\pi^m_n(A)>>\mathcal{L}_n(X)(A)
\end{CD}$$
commutes.

Since the schemes $\mathcal{L}_n(X)$ are affine for affine $X$,
and $\mathcal{L}_n(\cdot)$ takes open covers to open covers, the
morphisms $\pi^m_n$ are affine for any $k$-variety $X$. This
property guarantees that the projective limit
$$\mathcal{L}(X)=\lim_{\stackrel{\longleftarrow}{n}}\mathcal{L}_n(X)$$ exists in the category of $k$-schemes.
 The scheme $\mathcal{L}(X)$ is called the arc
scheme of $X$. It is not Noetherian, in general. It comes with
natural projection morphisms
$$\pi_n:\mathcal{L}(X)\rightarrow \mathcal{L}_n(X)$$
 For any field $F$ over $k$, we have a natural bijection
$$\mathcal{L}(X)(F)=X(F[[t]])$$ and the points of these sets are
called $F$-valued arcs on $X$. The morphism $\pi_n$ maps an arc to
its truncation modulo $t^{n+1}$ in $X(F[t]/(t^{n+1}))$. In
particular, $\pi_0$ sends an element $\psi$ of $X(F[[t]])$ to the
element $\psi(0)$ of $X(F)$ obtained by putting $t=0$ in the
coordinates of $\psi$. We call $\psi(0)$ the origin of the arc
$\psi$. An arc should be seen as a two-dimensional infinitesimal
disc on $X$ with origin at $\psi(0)$.

It follows immediately from the definition that
$\mathcal{L}_0(X)=X$ and that $\mathcal{L}_1(X)$ is the tangent
scheme of $X$.  A morphism of $k$-varieties $h:Y\rightarrow X$
induces morphisms
\begin{eqnarray*}
\mathcal{L}(h)&:&\mathcal{L}(Y)\rightarrow
\mathcal{L}(X)\\\mathcal{L}_n(h)&:&\mathcal{L}_n(Y)\rightarrow
\mathcal{L}_n(X)\end{eqnarray*} which commute with the truncation
maps.

 If $X$ is smooth over $k$, of pure dimension $d$,
then the morphisms $\pi^m_n$ are Zariski-locally trivial
fibrations with fiber $\A^{d(m-n)}_k$. To see this, note that by
smoothness, $X$ can be covered with open subvarieties $U$ which
admit an \'etale morphism $h:U\rightarrow \A^d_k$. Since
$A[t]/(t^{m+1})\rightarrow A[t]/(t^{n+1})$ is a nilpotent
immersion for any $k$-algebra $A$, the infinitesimal lifting
criterion for \'etale morphisms implies that the square
$$\begin{CD}
\mathcal{L}_m(U)@>\mathcal{L}_m(h)>> \mathcal{L}_m(\A^d_k)
\\@V\pi^m_nVV @VV\pi^m_nV
\\ \mathcal{L}_n(U)@>\mathcal{L}_n(h)>> \mathcal{L}_n(\A^d_k)
\end{CD}$$
is Cartesian. In intuitive terms, arcs are (\'etale-)local objects
on $X$ and any smooth variety of pure dimension $d$ looks
(\'etale-)locally like an open subvariety of $\A^d_k$; but an
element of
$$\mathcal{L}_n(\A^d_k)(A)=\A^d_k(A[t]/t^{n+1})$$ is simply a
$d$-tuple of elements in $A[t]/t^{n+1}$.

If $X$ is singular, the schemes $\mathcal{L}_n(X)$ and
$\mathcal{L}(X)$ are still quite mysterious. They contain a lot of
information on the singularities of $X$. Motivic integration
provides a powerful way to extract interesting invariants of the
singularities from the arc schemes $\mathcal{L}(X)$.

\begin{example} We continue our previous Example. For any $k$-algebra
$A$, a point on $\mathcal{L}(X)$ with coordinates in $A$ is given
by a couple
$$(x(t)=x_0+x_1t+x_2t^2+\ldots,y(t)=y_0+y_1t+y_2t^2+\ldots)$$
 with $x_i,\,y_i\in A$, such that $x(t)^2-y(t)^3=0$.

 Working this out yields an infinite number of polynomial equations in the variables
 $x_i,\,y_i$ and these realize $\mathcal{L}(X)$ as a closed subscheme of the infinite-dimensional affine space
 $$\A^{\infty}_k=\Spec k[x_0,y_0,x_1,y_1,\ldots]$$
 The truncation map
 $$\pi_n:\mathcal{L}(X)\rightarrow \mathcal{L}_n(X)$$ sends
 $(x(t),y(t))$ to
 $$(x_0+\ldots+x_nt^n,y_0+\ldots+y_nt^n)$$ and (if $A$ is a field)
 the origin of $(x(t),y(t))$ is simply the point $(x_0,y_0)$ in $X(A)$.
 \end{example}

 If $k$ has characteristic zero, one can give
 an elegant construction of the schemes $\mathcal{L}_n(X)$ and
 $\mathcal{L}(X)$ using differential algebra. Assume that $X$ is
 affine, say given by polynomial equations
 $$f_1(x_{1},\ldots,x_r)=\ldots=f_\ell(x_1,\ldots,x_r)$$ in affine
 $r$-space $\A^r_k=\Spec k[x_1,\ldots,x_r]$. Consider the
$k$-algebra $$B= k[y_{1,0},\ldots,y_{r,0},y_{1,1},\ldots]$$ and
the unique $k$-derivation $\delta:B\rightarrow B$ mapping
$y_{i,j}$ to $y_{i,j+1}$ for each $i,j$. Then $\mathcal{L}(X)$ is
isomorphic to the closed subscheme of $\Spec B$ defined by the
equations
$$\delta^{(i)}(f_q(y_{1,0},\ldots,y_{r,0}))=0$$ for $q=1,\ldots,\ell$
and $i\in \N$. The point with coordinates $y_{i,j}$ corresponds to
the arc
$$(\sum_{j\geq 0}\frac{y_{1,j}}{j!}t^j,\ldots,\sum_{j\geq 0}\frac{y_{r,j}}{j!}t^j )$$
See for instance \cite{NiSe-gk} for details.

\subsection{Motivic integrals}\label{subsec-motint}
Copying the notion of cylinder and the description of its Haar
measure to the setting of arc spaces, we can define a
\textit{motivic measure} on a class of subsets of the arc space
$\mathcal{L}(X)$. From now on, we assume that $X$ is smooth over
$k$, of pure dimension $d$.

\begin{definition} A cylinder
in $\mathcal{L}(X)$ is a subset $C$ of the form
$(\pi_m)^{-1}(C_m)$, with $m\geq 0$ and $C_m$ a constructible
subset of $\mathcal{L}_m(X)$.
\end{definition}

 Note that the set of cylinders in $\mathcal{L}(X)$ is a
Boolean algebra, i.e. it is closed under complements, finite
unions and finite intersections.

\begin{lemma-def} Let $C$ be a
cylinder in $\mathcal{L}(X)$, and choose $m\geq 0$ such that
$C=(\pi_{m})^{-1}(C_m)$ with $C_m$ constructible in
$\mathcal{L}_m(X)$. The value
$$\mu(C):=[\pi_m(C)]\LL^{-d(m+1)} \quad \in \mathcal{M}_k$$
does not depend on $m$, and is called the motivic measure $\mu(C)$
of $C$.
\end{lemma-def}
\begin{proof} This follows immediately from the fact that the truncation
morphisms $\pi^{n}_{m}$ are Zariski-locally trivial fibrations
with fiber $\A^{d(n-m)}_k$.\end{proof} Here is an elementary but
important example.

\begin{example} If
$C=\mathcal{L}(X)$, then $\mu(C)=\LL^{-d}[X]$. \end{example}

 The normalization factor $\LL^{-d}$ is added in
accordance with the $p$-adic case, where the ring of integers gets
measure one (rather then the cardinality of the residue field).
For geometric applications, it would be more natural to omit it
(and this is often done in literature).
 Note that the motivic measure $\mu$ is additive w.r.t.
finite disjoint unions.

\begin{remark}
 In the general theory of motivic integration on algebraic varieties \cite{DLinvent}, one
constructs a much larger class of measurable subsets of arc spaces
of (possibly singular) $k$-varieties, and one defines the motivic
measure via approximation by cylinders. This necessitates
replacing $\mathcal{M}_k$ by a certain ``dimensional completion''
$\widehat{\mathcal{M}}_k$.
\end{remark}

The following definition suggests itself.
\begin{definition} We say that a function
$$\alpha:\mathcal{L}(X)\rightarrow \N\cup\{\infty\}$$
is integrable if it takes finitely many values, and if the fiber
$\alpha^{-1}(i)$ is a cylinder for each $i\in \N$.

 We define the motivic integral of $\alpha$ by
$$\int_{\mathcal{L}(X)}\LL^{-\alpha}=\sum_{i\in
\N}\mu(\alpha^{-1}(i))\LL^{-i}\quad \in \mathcal{M}_k$$
\end{definition}

 We will need the following generalization of the
definition.
\begin{definition} Let $s$ be a
formal variable, and consider functions
$$\alpha,\,\beta:\mathcal{L}(X)\rightarrow \N\cup\{\infty\}$$
We say that $\alpha\cdot s+\beta$ is integrable if
$\alpha^{-1}(i)$ and $\beta^{-1}(i)$ are cylinders for every $i\in
\N$, and if $\beta$ takes only a finite number of values on each
fiber $\alpha^{-1}(i)$ with $i\in \N$.

 We define the motivic integral of $\alpha\cdot
s+\beta$ by
$$\int_{\mathcal{L}(X)}\LL^{-(\alpha\cdot s+\beta)}=\sum_{i,j\in \N}\mu(\alpha^{-1}(i)\cap \beta^{-1}(j))\LL^{-(is+j)}\quad \in
\mathcal{M}_k[[\LL^{-s}]]$$
\end{definition}
 The condition that $\beta$ assumes only finitely many
values on the fibers of $\alpha$ guarantees that the coefficient
of each $\LL^{-is}$ in the definition is a finite sum.

\subsection{Change of variables}\label{subsec-changevar}
The central and most powerful tool in the theory of motivic
integration is the change of variables formula. It has various
profound applications in birational geometry and singularity
theory. For its precise statement, we need some auxiliary
notation. For any $k$-variety $Y$, any ideal sheaf $\mathscr{J}$
on $Y$ and any arc
$$\psi:\Spec F[[t]]\rightarrow Y$$ on $Y$, we define the order of $\mathscr{J}$ at $\psi$ by
 $$ord_t\mathscr{J}(\psi)=\min\{ord_t \psi^*f\,|\,f\in \mathscr{J}_{\psi(0)}\}$$
 where $ord_t\psi^*f$ is the $t$-adic valuation of $\psi^*f\in
 F[[t]]$. In this way, we obtain a function
$$ord_t\mathscr{J}:\mathcal{L}(Y)\rightarrow \N\cup\{\infty\}$$
 whose fibers over $\N$ are cylinders. Note that
 $ord_t\mathscr{J}(\psi)=\infty$ iff the image of $\psi$ is contained
 in the support of $\mathscr{J}$.

\begin{theorem}[Denef-Loeser; Change of variables
formula]\label{thm-changevar} Assume that $k$ is perfect. Let
$h:Y\rightarrow X$ be a proper birational morphism, with $Y$
smooth over $k$, and denote by $Jac_h$ the Jacobian ideal sheaf of
$h$. Let $\alpha\cdot s+\beta$ be an integrable function on
$\mathcal{L}(X)$, and assume that $ord_tJac_h$ takes only finitely
many values on each fiber of $\alpha\circ \mathcal{L}(h)$ over
$\N$. Then
$$\int_{\mathcal{L}(X)}\LL^{-(\alpha\cdot s+\beta)}=\int_{\mathcal{L}(Y)}\LL^{-((\alpha\circ \mathcal{L}(h))\cdot s+
\beta\circ \mathcal{L}(h)+ord_tJac_h)}$$ in
$\mathcal{M}_k[[\LL^{-s}]]$.
 \end{theorem}
 \begin{remark}
In \cite{DLinvent}, Denef and Loeser stated the change of
variables formula in the case where $k$ has characteristic zero,
but this assumption is not necessary; see \cite{sebag1}.
 \end{remark}

The proof of the change of variables formula goes beyond the scope
of this introduction. The very basic idea behind the formula is
the following:
 if we denote by $V$ the closed subscheme of $Y$ defined by the Jacobian ideal $Jac_h$, and by $U$
 its image under $h$, then the morphism $h:Y-V\rightarrow X-U$ is an
 isomorphism. Combined with the properness of $h$, this implies that
$$\mathcal{L}(h):\mathcal{L}(Y)-\mathcal{L}(V)\rightarrow
\mathcal{L}(X)-\mathcal{L}(U)$$ is a bijection; but
$\mathcal{L}(V)$ and $\mathcal{L}(U)$ have measure zero in
$\mathcal{L}(Y)$, resp. $\mathcal{L}(X)$ (w.r.t. a certain more
refined motivic measure) so it is reasonable to expect that there
exists a change of variables formula.

 The jet schemes $\mathcal{L}_n(Y)$, however, are ``contracted'' under the morphism
 $$\mathcal{L}_n(h):\mathcal{L}_n(Y)\rightarrow \mathcal{L}_n(X)$$ and this affects the motivic measure of cylinders.
 The ``contraction factor'' is measured by the Jacobian.

 \begin{example}
Put $X=\A^3_k=\Spec k[x,y,z]$ and let $h:Y\rightarrow X$ be the
blow-up of $X$ at the origin $O$. Let $\mathscr{J}$ be the ideal
sheaf on $X$ defining the origin. It corresponds to the ideal
$(x,y,z)$ in $k[x,y,z]$. Denote by $C$ the subset of
$\mathcal{L}(X)$ where the function $ord_t\mathscr{J}$ takes the
value $1$. It consists of the arcs
$$(x_0+x_1t+x_2t^2+\ldots,y_0+y_1t+y_2t^2+\ldots,z_0+z_1t+z_2t^2+\ldots)$$ on $\A^3_k$ with
$x_0=y_0=z_0=0$ and $(x_1,y_1,z_1)\neq (0,0,0)$, so we see that it
is a cylinder over the constructible subset $\pi_1(C)$ of
$$\mathcal{L}_1(X)=\Spec k[x_0,x_1,y_0,y_1,z_0,z_1]$$ defined by
the same conditions. Since $\pi_1(C)$ is isomorphic to
$\A^3_k\setminus \{(0,0,0)\}$, the motivic measure of $C$ is equal
to
$$\mu(C)=\LL^{-6}[\pi_1(C)]=\LL^{-6}(\LL^3-1)$$ in
$\mathcal{M}_k$.

The pull-back $\mathscr{J}'=\mathscr{J}\mathcal{O}_Y$ of the ideal
sheaf $\mathscr{J}$ to $Y$ is the defining ideal sheaf of the
exceptional divisor $E=h^{-1}(O)$. The jacobian ideal $Jac_h$ of
$h$ is the square of $\mathscr{J}'$. The inverse image $C'$ of $C$
under the morphism $\mathcal{L}(h)$ is the cylinder in
$\mathcal{L}(Y)$ where the function $$ord_t
\mathscr{J}'=ord_t\mathscr{J}\circ\mathcal{L}(h)$$ takes value
$1$. This implies that $ord_t Jac_h$ takes value $2$ on $C'$.

 The set $C'$ is a cylinder over its
image $\pi_1(C')\cong E\times \A^2_k\times_k\mathbb{G}_{m,k}$ in
$\mathcal{L}_1(C')$. Indeed: if $u$ is any closed point of $E$ and
we choose local coordinates $(x',y',z')$ on an open neighbourhood
$U$ of $u$ (i.e. an \'etale morphism $U\rightarrow \A^3_k=\Spec
k[x',y',z']$ mapping $u$ to the origin) such that $E\cap U$ is
defined by the equation $x'=0$, then $C'\cap (\pi_0)^{-1}(U)$ is
determined by the conditions $x_0'=0$ and $x_1'\neq 0$.

The motivic measure of $C'$ is equal to
$$\mu(C')=\LL^{-6}(\LL-1)\LL^2[E]$$ in $\mathcal{M}_k$. The change of variables formula tells us that
$$\mu(C)=\mu(C')\LL^{-2}$$ which is indeed the case since $E$ is
isomorphic to $\Pro^2_k$ and $[E]=\LL^2+\LL +1$ in
$\mathcal{M}_k$.
 \end{example}

 \subsection{A (cheating) proof of Kontsevich' theorem}
As we already mentioned, Kontsevich created motivic integration to
prove the following theorem. Recall that a Calabi-Yau variety over
$\C$ is a smooth, proper, connected $\C$-variety $X$ of dimension
$d$ such that $\Omega^d_{X/\C}$ is isomorphic to $\mathcal{O}_X$,
i.e. such that $X$ admits a nowhere vanishing differential form of
maximal degree.
\begin{theorem}[Kontevich]
If $X$ and $Y$ are birationally equivalent Calabi-Yau varieties
over $\C$, then $X$ and $Y$ have the same Hodge numbers.
\end{theorem}

To prove this result, we need the refined version of motivic
integration taking values in a certain completion
$\widehat{\mathcal{M}}_{\C}$ of $\mathcal{M}_{\C}$. Nevertheless,
it is possible to indicate at least the general ideas. Denote by
$d$ the dimension of $X$ and $Y$. The fact that $X$ and $Y$ are
birationally equivalent implies that there exists a smooth and
proper $\C$-variety $Z$ together with proper birational morphisms
$$\begin{CD}X@<f<< Z@>g>> Y\end{CD}$$
We can express $\LL^{-d}[X]\in \mathcal{M}_{\C}$ by the motivic
integral
$$\int_{\mathcal{L}(X)}1$$
(where we should really write $\LL^{-0}$ instead of $1$) and the
analogous expression holds for $\LL^{-d}[Y]$. Now we can compute
both motivic integrals on $Z$, using the change of variables
formula. This yields
\begin{eqnarray}
\LL^{-d}[X]&=&\int_{\mathcal{L}(Z)}\LL^{-ord_tJac_f} \label{X}
\\ \LL^{-d}[Y]&=& \int_{\mathcal{L}(Z)}\LL^{-ord_tJac_g} \label{Y}
\end{eqnarray}
The alert reader may have noticed that this is the place where
we're cheating, since $ord_tJac_f$ and $ord_tJac_g$ take
infinitely many values on $\mathcal{L}(Z)$ if $f$ and $g$ are not
isomorphisms, so that the integrals on the right hand side become
infinite sums. However, the general formalism of motivic
integration can deal with such expressions, replacing
$\mathcal{M}_{\C}$ by $\widehat{\mathcal{M}}_{\C}$.

One can show that the specialization morphism $HD$ from Section
\ref{subsec-groth} factors through the image of $\mathcal{M}_{\C}$
in $\widehat{\mathcal{M}}_{\C}$, so to prove the theorem it
suffices to show that (\ref{X}) and (\ref{Y}) coincide. Choose
nowhere vanishing differential forms $\omega_X$ and $\omega_Y$ of
degree $d$ on $X$, resp. $Y$. Since $X$ is normal and $Y$ is
proper, the rational map $g\circ f^{-1}$ is defined on an open
subvariety of $X$ whose complement has codimension $\geq 2$ in
$X$, and the pullback $(g\circ f^{-1})^*\omega_Y$ extends uniquely
to a degree $d$ differential form on $X$. Hence, there exists a
unique regular function $a$ on $X$ with $(g\circ
f^{-1})^*\omega_Y=a\cdot \omega_X$. Pulling back to $Z$ we get
$g^*\omega_Y=f^*a\cdot f^*\omega_X$. Symmetrically, there exists a
unique regular function $b$ on $Y$ such that $(f\circ
g^{-1})^*\omega_X=b\cdot \omega_Y$ and on $Z$ we have
$f^*\omega_X=g^*b\cdot g^*\omega_Y$. Hence, $f^*a$ and $g^*b$ are
units on $Z$ and belong to $\C$. This means that the ideal sheaves
$Jac_f$ and $Jac_g$ coincide.


A similar proof yields a stronger result: we say that two proper,
smooth, connected $\C$-varieties $X$ and $Y$ are $K$-equivalent if
there exists a smooth, proper, connected $\C$-variety $Z$ and
proper birational morphisms $f:Z\rightarrow X$ and $g:Z\rightarrow
Y$ such that the Jacobian divisors $(Jac_f)$ and $(Jac_g)$ are
linearly equivalent (i.e. $Jac_f$ and $Jac_g$ are isomorphic as
invertible sheaves on $Z$). One can show that this automatically
implies that the ideal sheaves $Jac_f$ and $Jac_g$ coincide, and
the above arguments show that $[X]=[Y]$ in
$\widehat{\mathcal{M}}_{\C}$.
\section{Motivic zeta functions}\label{sec-motzeta}
\subsection{Definitions}\label{subsec-defmotzeta}
 Using the
ideas in the previous section, we can transfer the definition of
the $p$-adic local zeta function to the motivic framework. We will
only consider the case where $\Phi$ is the characteristic function
of $R^d$, and the character $\chi$ is trivial, i.e. we'll
construct the motivic counterpart of
$$Z(f;s)=\int_{R^d}|f|_K^sd\mu$$
with $f\in K[\widetilde{x}]=K[x_1,\ldots,x_d]$.
  There is also a more delicate definition of the motivic
zeta function for non-trivial characters $\chi$ and more general
$\Phi$; see \cite{DL5}.

Let $k$ be any field, and let $f$ be a polynomial in
$k[x_1,\ldots,x_d]$. The integration space $R^d$ is replaced by
the arc scheme $\mathcal{L}(\A^d_k)$ with the motivic measure.
What should $|f|_K^s$ be?

For $g\in K[\widetilde{x}]$ and $x\in R^d$, $|g(x)|_K$ equals
$q^{-v(g(x))}$, with $q$ the number of rational points on the
affine line over the residue field of $R$.
 So in the motivic setting, the natural candidate to replace
$|f(x)|_K$ for some $k$-field $F$ and $$x\in
\mathcal{L}(\A^d_k)(F)=F[[t]]^d$$
 is
$\LL^{-ord_t f(x)}$ where $ord_t$ denotes the $t$-adic valuation
on  $F[[t]]$. This leads to the following definition.

\begin{definition} Let $k$ be any field, let $X$ be a
smooth $k$-variety of pure dimension, and let
$$f:X\rightarrow \A^1_k$$ be a $k$-morphism to the affine line. We
define the motivic zeta function $Z_{mot}(f;s)$ of $f$ by
$$Z_{mot}(f;s)=\int_{\mathcal{L}(X)}\LL^{-ord_tf\cdot s}\quad \in
\mathcal{M}_k[[\LL^{-s}]]$$
\end{definition}
Here  $ord_tf=ord_t(f)$, the function $\mathcal{L}(X)\rightarrow
\N\cup\{\infty\}$ associated to the ideal sheaf $(f)$ on $X$
generated by $f$.

Denote by $V_f$ the closed subscheme of $X$ defined by $(f)$. As
in the $p$-adic case, a simple transformation rule relates this
zeta function to the motivic Poincar\'e series
$$Q_{mot}(T)=\sum_{m\geq 0}[\mathcal{L}_m(V_f)]T^{m+1}\quad \in K_0(Var_k)[[T]]$$
with $T=\LL^{-s}$. This is the motivic counterpart of Igusa's
Poincar\'e series, since for any $k$-field $F$, the set
$\mathcal{L}_m(V_f)(F)$ is the space of solutions of the
congruence $f\equiv 0\mod t^{n+1}$ in $X(F[t]/(t^{n+1}))$.

\begin{theorem}[Denef-Loeser]
Suppose that $k$ has characteristic zero. The motivic zeta
function $Z_{mot}(f;s)$ is rational in $T=\LL^{-s}$ over
$\mathcal{M}_k$. More precisely, there exists a finite subset $S$
of $\Z_{>0}\times \Z_{>0}$ such that $Z_{mot}(f;s)$ belongs to
$$\mathcal{M}_k\left[\frac{\LL^{-as-b}}{1-\LL^{-as-b}} \right]_{(a,b)\in S} \subset
\mathcal{M}_k[[\LL^{-s}]]$$
\end{theorem}
 If $h:Y\rightarrow X$ is an embedded resolution of singularities for $f$,
then we can take for $S$ the set of numerical data of the
resolution. In fact, by a similar computation as in the $p$-adic
case, using the change of variables formula, one obtains an
explicit expression for $Z_{mot}(f;s)$ in terms of $h$. Note that
Theorem \ref{thm-changevar} applies to the morphism $h$~: since we
always assume that $h$ is an isomorphism over the complement of
the zero locus of $f$ in $X$, the Jacobian divisor $(Jac_h)$ is
supported on $(f\circ h)$, which implies that $ord_tJac_h$ takes
only finitely many values on each fiber of $ord_t(f\circ h)=ord_t
f \circ \mathcal{L}(h)$.

Denote by $E_i,\,i\in I$ the irreducible components of the divisor
$E=(f\circ h)$, and by $(N_i,\nu_i)$ the corresponding numerical
data. For any non-empty subset $J$ of $I$ we put
 $E_J=\cap_{j\in J}E_j$ and $E_J^o=E_J\setminus(\cup_{i\notin
 J}E_i)$. In particular,  $E_{\emptyset}=Y$ and $E_{\emptyset}^o=Y\setminus
 E$. We denote by $d$ the dimension of $X$.

 \begin{theorem}[Denef-Loeser]\label{theo-explicit}
$$Z_{mot}(f;s)=\LL^{-d}\sum_{J\subset
I}[E_J^o]\prod_{j\in
J}\frac{(\LL-1)\LL^{-N_is-\nu_i}}{1-\LL^{-N_is-\nu_i}}\quad \in
\mathcal{M}_k[[\LL^{-s}]]$$
 \end{theorem}
\subsection{The motivic monodromy
conjecture}\label{subsec-motmoncon}
 If we want to
translate the statement of the monodromy conjecture to the motivic
setting, there is a technical complication: one should be careful
when speaking about poles of the zeta function, since $K_0(Var_k)$
is not a domain.
\begin{conjecture}[Denef-Loeser; Motivic monodromy conjecture]
Suppose that $k$ is a subfield of $\C$. Let $X$ be a smooth
$k$-variety and $f:X\rightarrow \A^1_k$ a morphism of
$k$-varieties, and denote by $X_s$ the complex analytic space
defined by the equation $f=0$ on the analytic space $X^{an}$
associated to the complex variety $X\times_k \C$.

There exists a finite subset $S$ of $\Z_{>0}\times \Z_{>0}$ such
that
$$Z_{mot}(f;s)\in \mathcal{M}_k\left[\LL^{-s},\frac{1}{1-\LL^{-as-b}}
\right]_{(a,b)\in S}$$
 and such that for each $(a,b)\in S$, the
value $\exp(-2\pi i b/a)$ is a monodromy eigenvalue of $f$ at some
point of $X_s$.
\end{conjecture}
 We will see that the motivic monodromy conjecture implies
the $p$-adic one.
\subsection{Specialization to the $p$-adic
world}\label{subsec-special} We'll now explain how the motivic
zeta function relates to the $p$-adic one. Let $L$ be a number
field, and let $f$ be an element of $L[x_1,\ldots,x_d]$. The
polynomial $f$ defines a morphism of $L$-varieties
$f:\A^d_L\rightarrow \A^1_L$. We'll denote by $\mathcal{O}_{L}$
the ring of integers of $L$, and by $\mathrm{Msp}\,\mathcal{O}_L$
its maximal spectrum. For each finite place $\mathfrak{P}\in
\mathrm{Msp}\,\mathcal{O}_L$, we denote its residue field by
$k_{\mathfrak{P}}$.

 On the one hand, we can associate to $f$ its motivic zeta
function $$Z_{mot}(f;s)\in \mathcal{M}_L[[\LL^{-s}]]$$ On the
other hand, for each finite place $\mathfrak{P}$, we can consider
the $\mathfrak{P}$-adic zeta function $Z_{\mathfrak{P}}(f;s)$.

Consider the ring
$$\mathscr{Z}=\left(\prod_{\mathfrak{P}\in \mathrm{Msp\,\mathcal{O}_L}} \Q \right)/
\left(\bigoplus_{\mathfrak{P}\in \mathrm{Msp\,\mathcal{O}_L}} \Q
\right)$$ For any variety $X$ over $L$, we can choose a model over
 $\mathcal{O}_L$, and count rational points on
the reduction modulo $\mathfrak{P}$, for each finite place
$\mathfrak{P}$. The outcome may depend on the chosen model, but
all these values together yield a well-defined element of
$\mathscr{Z}$. Moreover, since this operation is additive (w.r.t.
closed immersions) and multiplicative (w.r.t. the fiber product
over $L$), we obtain a morphism of rings $K_0(Var_L)\rightarrow
\mathscr{Z}$ which induces a morphism of rings
$$\mathscr{N}:\mathcal{M}_k\left[\frac{\LL^{-as-b}}{1-\LL^{-as-b}} \right]_{(a,b)\in \Z_{>0}\times \Z_{>0}}
\rightarrow \mathscr{Z}'$$ with
$$\mathscr{Z}'=\left(\prod_{\mathfrak{P}\in \mathrm{Msp\,\mathcal{O}_L}} \Q\left[\frac{|k_{\mathfrak{P}}|^{-as-b}}{1-|k_{\mathfrak{P}}|^{-as-b}} \right]_{(a,b)\in
\Z_{>0}\times \Z_{>0}} \right)/ (\oplus \ldots)$$

\begin{theorem}[Denef-Loeser]\label{specialize} The
motivic zeta function $Z_{mot}(f;s)$ specializes to
$Z_{\mathfrak{P}}(f;s)$ for almost all finite places
$\mathfrak{P}$, in the following sense: the image of
$Z_{mot}(f;s)$ under the morphism $\mathscr{N}$ is the quotient
class of the tuple
$$(Z_{\mathfrak{P}}(f;s))_{\mathfrak{P}\in
\mathrm{Msp}\,\mathcal{O}_L}$$
\end{theorem}
\begin{proof}
Combine the expressions in Theorem \ref{theo-denef} and Theorem
\ref{theo-explicit}.
\end{proof}

 \begin{cor} The motivic monodromy conjecture for $k=L$ and
$f\in L[\widetilde{x}]$ implies the untwisted $p$-adic monodromy
conjecture for $f$ (Conjecture \ref{con-untwist}).
\end{cor}

In fact, Denef and Loeser formulated a more general motivic
monodromy conjecture, involving motivic zeta functions which are
twisted by characters $\chi$. A similar specialization result
holds in that setting; see \cite{DL5}.

 In virtually all of the cases where the $p$-adic monodromy
conjecture is proven, the same strategy yields a proof for the
motivic version of the conjecture.
\subsection{Why motivic zeta functions?}\label{subsec-why}
The motivic zeta function appears in Theorem \ref{specialize} as a
universal zeta function, with the $p$-adic zeta functions as its
avatars. It captures the geometric structure which explains the
uniform behaviour of the $p$-adic zeta functions for $p\gg 0$. In
this sense it fully deserves the name ``motivic''.

Although more general than the $p$-adic one, the motivic mondromy
conjecture seems more accessible, since we never leave the
equicharacteristic zero world. The arc spaces appearing in its
definition are closely related to the infinitesimal structure of
the morphism $f$, and hence to its singularities, so one can
certainly believe that something like the motivic monodromy
conjecture should hold. The connection becomes even stronger when
we consider the so-called monodromic motivic zeta function
\cite[3.2]{DL3}, where we actually see the monodromy appear in the
form of an action of the profinite group $\hat{\mu}$ of roots of
unity on the coefficients of the zeta function. Intriguingly, by
taking a formal limit of this monodromic zeta function one gets a
motivic incarnation of the nearby cycles of $f$, and the monodromy
eigenvalues can be read from this object \cite[3.5]{DL3}.
Nevertheless, at this moment the conjecture still seems far out of
reach.

The motivic zeta functions are also interesting in their own
right, independently of the relation with $p$-adic zeta functions
and the monodromy conjecture. They provide very fine invariants of
hypersurface singularities, which can be explicitly computed on a
resolution of singularities. On the other hand, since motivic zeta
functions are defined intrinsically by means of a motivic
integral, they show that certain invariants of a resolution of
singularities are actually independent of the chosen resolution.
See for instance \cite[\S\,1]{loeser-seattle} for the genesis of
these ideas.

\end{document}